\documentclass[oneside,reqno, twoside, 11pt]{amsart}
\usepackage{color}
\usepackage[T1]{fontenc}
\usepackage{lmodern}
\usepackage[english]{babel}
\usepackage{mathrsfs}  
\usepackage{upgreek}
\usepackage{amsfonts}
\usepackage{extpfeil}
\usepackage{verbatim}
\usepackage{BOONDOX-cal}
\usepackage{microtype}
\usepackage[new]{old-arrows}
\usepackage{amsmath}
\usepackage[all, cmtip]{xy}
\usepackage{tikz-cd}

\usepackage{amsthm}
\usepackage{mathtools}

\usepackage{bbm}
\usepackage{paralist}
\usepackage[T1]{fontenc}
\usepackage[colorlinks,citecolor=magenta,pagebackref=true,urlcolor=magenta, bookmarksdepth=3]{hyperref}
\usepackage[nodisplayskipstretch]{setspace}
\usepackage{breqn}

\usepackage{nicefrac}
\usepackage[font=small,labelfont=bf]{caption}
\usepackage[labelformat=simple]{subcaption}

\usepackage{epsfig}

\DeclareFontFamily{OT1}{pzc}{}
\DeclareFontShape{OT1}{pzc}{m}{it}{<-> s * [1.10] pzcmi7t}{}
\DeclareMathAlphabet{\mathpzc}{OT1}{pzc}{m}{it}

\makeatletter
\newtheorem*{rep@theorem}{\rep@title}
\newcommand{\newreptheorem}[2]{%
	\newenvironment{rep#1}[1]{%
		\def\rep@title{#2~\ref{##1}}%
		\begin{rep@theorem}}%
		{\end{rep@theorem}}}

\makeatother
\theoremstyle{plain}
\newtheorem*{thm*}{Theorem}
\newtheorem{thm}{Theorem}[section]
\newreptheorem{mthm}{Main Theorem}
\newreptheorem{mcor}{Corollary}
\newreptheorem{mlem}{Lemma}
\newreptheorem{thm}{Theorem}

\newtheorem{lem}[thm]{Lemma}

\newtheorem*{lem*}{Lemma}

\theoremstyle{definition}

\newtheorem{ex}[thm]{Example}

\newcommand{\supp}{\mathrm{supp}\hspace{1mm}}

 \newcommand{\R}{\mathbb{R}}

\newcommand{\tet}{\theta}

\newcommand{\st}{\mr{st}}

\newcommand{\sbseq}{\subseteq}
\newcommand{\spseq}{\supseteq}

\newcommand{\vanish}[1]{}

\def\V{{\bf V}}

\def\im{\mathrm{im}}%\hspace{1mm}}
\def\ker{\mathrm{ker}}%\hspace{1mm}}

\def\sbs\subset
\def\sbseq{\subseteq}

\def\langle{\left<}
\def\rangle{\right>}

\def\({\left(}
\def\){\right)}
\def\no={\,{\,|\!\!\!\!\!=\,\,}}

\def\no={\,{\,|\!\!\!\!\!=\,\,}}
\def\sbseq{\subseteq}

\def\sbseq{\subseteq}
\def\sbs\subset
\def\spseq{\supseteq}

\newcommand{\xqedhere}[2]{%
	\rlap{\hbox to#1{\hfil\llap{\ensuremath{#2}}}}}

\newcommand\Defn[1]{\textbf{#1}}
\newcommand{\cm}[1]{}

\newcommand\mbf[1]{\mathbf{#1}}

\newcommand\mr[1]{\mathrm{#1}}

\newcommand\x{\mathbf{x}}

\DeclareMathOperator{\Lk}{lk}

\DeclareMathOperator{\St}{st}

\title{Unstable blueprints can be shared}
\author{Karim Adiprasito }

\date{24.06.2019}

%Primary  14M25, 05C38 ; Secondary 32S50, 52C25,  13F55}

\begin{document}
	
	\begin{abstract}
This expository note illustrates toric perturbation and biased pairing theory to show that Artinian reductions of face rings of $2$-spheres that do not satisfy the Lefschetz property can be cut along a flat equator. This complements classical work of Bricard and Connelly, and exhibits a fundamental symmetry in non-rigid triangulations of spheres.
	\end{abstract}
	
	\maketitle
	
	\newcommand{\AR}{\mathcal{A}}
	\newcommand{\BR}{\mathcal{B}}
	\newcommand{\CR}{\mathcal{C}}
	\newcommand{\Mu}{M}
	\newcommand{\Soc}{\mathcal{S}\hspace{-1mm}\mathcal{o}\hspace{-1mm}\mathcal{c}}
	\newcommand{\Socl}{{\Soc^\circ}}

\section{Rigidity}

Cauchy proved that all triangulated $2$-spheres are rigid when they are realized as boundaries of 2-dimensional polytopes, that is, there is up to global isometry of $\R^3$ no other realization of that $2$-sphere as a convex $3$-polytope with the same edge-lengths.

Alexandrov showed then that triangulated $2$-spheres are infinitesimally rigid: He showed that there is no infinitesimal deformation of a simplicial $3$-polytope that does not extend to a global isometry or changes the lengths of edges in the first order. An immediate corollary is that triangulated $2$-spheres are generically rigid: a generic realization (choice of coordinates of the vertices) of the sphere in  $\R^3$ is infinitesimally rigid. 

On the other hand, infinitesimally flexible, and even flexible, realizations of $2$-dimensional triangulated spheres do exist, as constructed by Bricard and Connelly. We refer to \cite{Connelly:RigiditySurvey} for a survey of the facts presented here.

Another related question that was open until recently is the following variation on the existence of flexible polyhedra: Is there a triangulated $2$-sphere with coordinates in $\R^3$ such that no linear movement of the vertices along their linear spans generates a infinitesimally rigid sphere. In other words, if we think of the $2$-sphere as an architectural blueprint under the radial projection to $S^2$, then every lift of that blueprint to a building in $\R^3$ is infinitesimally rigid. To motivate this question, which we shall call the problem of unstable blueprints, in another way, it is useful to translate the problem to an algebraic setting.

\section{Rings}
If $\Delta$ is an abstract simplicial complex defined on the groundset $[n]\coloneqq \{1,\cdots,n\}$, let $I_\Delta\coloneqq \langle \x^{\mbf{a}}:\ \supp(\mbf{a})\notin\Delta\rangle$ denote the nonface ideal in $\R[\x]$, where $\R[\x]=\R[x_1,\cdots,x_n]$. Let $\R^\ast[\Delta]\coloneqq \R[\x]/I_\Delta$ denote the \Defn{face ring} of $\Delta$. Now, we pick a sufficient number of linear forms to make sure the quotient is finite dimensional:

We may associate to the vertices of $\Delta$ the coordinates $\V_\Delta=({v_1},\cdots, {v_n}) \in \R^{l\times n}$, obtaining a system of linear forms by considering  $\V_\Delta\x = \Theta$. With this, we obtain a \Defn{geometric simplicial complex}.

The face ring of a {geometric simplicial complex} $\Delta$ is considered with respect to its natural system of parameters induced by the coordinates, that is, 
\[\AR^\ast(\Delta)\ :=\ \R^\ast[\Delta]/\Theta\R^\ast[\Delta].\]

A geometric simplicial complex in $\R^d$ is \Defn{proper} if the image of every $k$-face, with $k<d$, linearly spans a subspace of dimension $k+1$. If $\Delta$ is of dimension $(d-1)$, and is given a proper coordinates in $\R^d$, then $\AR^\ast(\Delta)$ is finite-dimensional as a vector space. This is also called the Artinian reduction of the face ring $\R^\ast[\Delta]$.

The Lefschetz property for a properly realized triangulated sphere $\Sigma$ is then the existence of  $\ell$ in $\AR^1[\Sigma]$ so that
\[\AR^k(\Sigma)\ \xrightarrow{\ \cdot \ell^{d-2k} \ }\ \AR^{d-k}(\Sigma). \]
is an isomorphism for every $k\le \frac{d}{2}$. 

Following \cite{Lee}, the problem of unstable blueprints then asks for a proper realization of a $2$-sphere in $\R^3$ that do not have the Lefschetz property with respect to $k=1$. The existence of such spheres was unknown until recently \cite[Section 4.5]{AHL}, and is perhaps slightly surprising, given that the associated face rings are still Gorenstein. We show the following:

\begin{thm}
Every unstable blueprint $\Sigma$ contains a dividing equator, that is, an embedded simple cycle in its $1$-skeleton that lies in a linear hyperplane of $\R^3$.
\end{thm}

Unfortunately, this criterion is not sufficient: there are spheres with the Lefschetz property that do have dividing equator, such as the octahedron.

 We should note also that the theorem extends to surfaces, but this requires the methods of \cite{AHL}. To keep the argument self-contained, we do not go into that case here. 
 
 One could equally ask for good necessary or even sufficient conditions for Artinian reductions of general, higher-dimensional triangulations of spheres. Unfortunately, these are not so easy, and depend on pairing properties of monomial ideals with respect to the Poincar\'e pairing, see \cite{AJS}.
 
Only in dimension $2$ and lower does the problem seem to have a nice answer in terms of necessity. An intermediate problem that seems tractable is to find simple sufficient conditions for Artinian reductions of face rings of $2$-spheres to have the Lefschetz property.

\section{Notation} Recall that the  \Defn{star} and \Defn{link} of a face $\sigma$ in $\Delta$ are
the subcomplexes \[\St_\sigma \Delta\ \coloneqq \ \{\tau:\exists \tau'\supset \tau,\ \sigma\subset
\tau'\in \Delta\}\ \
\text{and}\ \ \Lk_\sigma \Delta\ \coloneqq \ \{\tau\setminus \sigma: \sigma\subset
\tau\in \Delta\}.\]
For geometric simplicial complexes $\Delta$, we shall think of the star of a face as a geometric subcomplex of $\Delta$, and the link of a face $\sigma$ as the geometric simplicial complex obtained by the orthogonal projection to $\mr{span}(\sigma)^\bot.$
Let us denote the \Defn{deletion} of $\sigma$ by $\Delta-\sigma$, the maximal subcomplex of $\Delta$ that does not contain $\sigma$.
Let 
\[\St_\sigma^\circ \Delta\ \coloneqq\  (\St_\sigma \Delta,\St_\sigma \Delta-\sigma).\]
We also recall that, for any vertex $v\in \Delta$, where $\Delta$ is a geometric simplicial complex in $\R^d$, and for any integer $k$, we have isomorphisms
	\[\AR^k({\Lk}_v \Delta)\ \cong\ \AR^{k}(\St_v \Delta)\ \xrightarrow{\ \cdot x_v \ }\ \AR^{k+1}(\St_v^\circ \Delta).\]

\section{A bad Artinian reduction} Before we continue, let us recall that the example of Artinian reduction without the Lefschetz property of \cite{AHL}: If $\Delta$ is a simplicial complex, and $\sigma$ is a face in $\Delta$ then a \Defn{stellar subdivision} of $\Delta$ at $\sigma$ is the simplicial complex
	\[\Delta{\uparrow}\sigma \ \coloneqq \ (\Delta-\sigma)\ \cup\ \bigcup_{\tau \in \partial \st_\sigma \Delta} (\{v_\sigma\}\ast \tau)
	\]
	where $v_\sigma$ is the new vertex introduced in place of $\sigma$, and $\ast$ is used to denote the free join of two simplices. If $\Delta$ is geometric, then $v_\sigma$ will be chosen to lie in the linear span of $\sigma$.
	
	\begin{ex}
		Consider $\Sigma$ the boundary of the tetrahedron, and denote its  vertices by $1,2,3$ and $4$.
		
		Perform stellar subdivisions at every triangle of $\Sigma$, call the newly created vertices $1',2',3'$ and $4'$. Call the resulting complex $\Sigma'$.
		
		Realize the vertices in $\R^{[3]}$ as follows: place the vertices $1,2,3$ and $4$ in general position in~$\R^{[2]}\subset \R^{[3]}$. Place the remaining vertices in $\R^{[3]}{\setminus} \R^{[2]}$. The associated linear system is a linear system of parameters for $\R[\Sigma']$, and therefore $\AR^\ast(\Sigma')$ is a Poincar\'e duality algebra.
		
		However, if we consider the subcomplex $\Delta=\Sigma \cap \Sigma'$,
		the quotient $\AR^\ast(\Delta)$ of $\AR^{\ast}(\Sigma')$ satisfies
		\[\dim \AR^1(\Delta)\ =\ 2\ \quad \text{and}\ \quad \dim \AR^2(\Delta)=3.\]
		But if $\Sigma$ has the Lefschetz property then
		\[\begin{tikzcd}
 \AR^{1}(\Sigma') \arrow{r}{\ \cdot \ell\ } \arrow[two heads]{d}{} & \AR^{2}(\Sigma') \arrow[two heads]{d}{} \\
 \AR^{1}(\Delta) \arrow{r}{\ \cdot \ell\ } & \AR^{2}(\Delta)
\end{tikzcd}
\]
The top horizontal map is an isomorphism, in particular a surjection. Hence, so is the bottom map, which contradicts the dimensions we computed for them.
		\end{ex}
		
		Let us briefly note that this almost trivial example provides also considerable strengthening of results of \cite{BH, ABA}: These papers provided minimal positively balanced fans that do not satisfy the Hodge-Riemann relations. By choosing all vertices of the stellar subdivisions to lie in one halfspace, we obtain positive balancing also here, but to a much stronger effect: not only does the Hodge-Riemann property fail, but even the Hard Lefschetz property does.

\section{Generic combinations of linear maps}

For the proof of the Main Theorem, we shall suppose the contrary. Assuming no dividing equator exists, we show the existence of a Lefschetz element. The idea to construct the map $\ell$ iteratively. We rely on the following principle:

\begin{lem}[Lemma 6.1, \cite{AHL}]\label{lem:perturbation}
	Given two linear maps 
	\[A, B: \mathcal{X}\ \longrightarrow\ \mathcal{Y}\]
	of two vector spaces $\mathcal{X}$ and $\mathcal{Y}$ over $\R$.
	Assume that 
		\[B(\ker A)\ \cap\ \im A\ =\ {0}\ \subset\ \mathcal{Y} .\]
		Then a generic linear combination 
		$A ``{+}" B$ of $A$ and $B$
		has kernel 
		\[\ker (
		A\ ``{+}"\ B)\ = \ \ker A\ \cap\ \ker B.\]
\end{lem}

The second observation is as follows:

\begin{lem}
A linear form $\ell$ on a zero-dimensional sphere satisfies the Lefschetz theorem if and only if $\ell$ is not linear.
\end{lem}

\section{Decomposing spheres}
Let us explore how these two lemmata provide Lefschetz elements, following \cite[Section 6.5]{AHL}: Assume we are in a sphere $\Sigma$ of dimension $d-1=2k$, and have shown that 

\begin{compactenum}[(1)]
\item There exists a set of vertices $W$ such that $\Sigma-W$ is a disk $\Delta$ and
\item We have the \Defn{transversal prime property} for $W$ in $\Sigma$, that is, the kernel of $\ell' = ``{\sum_{v\in W}}"x_v$ is exactly $\AR^k(\Delta,\partial \Delta),$ the intersection of the kernels of the $x_v, v\in W$.
\end{compactenum}

Lets pick $w$ a vertex in the boundary of $\Delta$. We wish to establish that $\ell'``+" x_w$ has kernel $\AR^k(\Delta-w,\partial (\Delta-w))$. 
Then we wish to understand what to do to ensure that the Poincar\'e pairing of $\AR^k(\st_w\Sigma)$ does not degenerate when restricting to the pullback of $\AR^k(\Delta,\partial \Delta)$. On the other hand,
we have a short exact sequence
\[\AR^k(\st_w \Delta, \st_w \partial \Delta)\ \xrightarrow{\ \cdot x_w\ }\ \AR^{k+1}(\Delta, \Delta-w)\ \longrightarrow\ \AR^{k+1}(\partial \Delta, \partial\Delta-w)\ \longrightarrow\ 0.\]
Hence, the kernel of the multiplication with $x_w$ is
$\AR^k(\Delta-w, \partial \Delta-w)$ provided \[\AR^{k+1}(\partial \Delta, \partial\Delta-w)\ =\ \AR^k({\Lk}_w \partial \Delta)\ =\ 0.\] Where is the Lefschetz theorem though? 

Well
$\uppi$ is a projection of $\R^d/w$ to a hyperplane, and $h$ is the coordinate with respect to that projection, $\tet=h\cdot\x\, $ the associated linear form, then the last condition is equivalent to the Lefschetz property in degree $k-1$: For the map
\[\AR^{k-1}(\uppi{\Lk}_w \partial \Delta)\ \xrightarrow{\ \cdot \tet \ }\ \AR^{k}(\uppi{\Lk}_w \partial \Delta)\]
to be an isomorphism is exactly equivalent to $\AR^k({\Lk}_w \partial \Delta)=0$.

\section{Proof of the Main Theorem}

Consider now a $2$-sphere $\Sigma$ in $\R^3$. We can start with any vertex to construct the Lefschetz element, and then have to iterate, finding iteratively new vertices such that

\begin{compactenum}[(1)]
\item The deletion of the initial segments of $i$ vertices from $\Sigma$ leaves $2$-disk $\Delta_i$, or is empty. This is to make sure we can keep inducting.
\item $\st_{i+1} \partial \Delta_i$ does not lie in a $2$-dimensional hyperplane. This is to ensure the transversal prime property when adding the vertex $i+1$, as follows from the previous two lemmas.
\end{compactenum}

We find new vertices iteratively. Assume now there is no dividing equator, then in particular, the boundary of  $\Delta_i$ does not lie in a linear hyperplane. In particular there exists a vertex $i+1$ satisfying (2) in its boundary. But that vertex, removed from $\\Delta_i$, might not leave a $2$-disk or be empty, but several disks joined along vertices.

Now, none of those smaller disks has a dividing equator boundary. We may therefore pick one of their boundary vertices instead, and choose one that satisfies (2). If those do not satisfy (1), we divided that disk into even smaller disks, and so on. Eventually, this process terminates, and we have found a way to find a new vertex that we can add to the already constructed set of vertices.

	\textbf{Acknowledgements:} Karim Adiprasito was supported by the European Research Council under the European
Union's Seventh Framework Programme ERC Grant agreement ERC StG 716424 - CASe, the Israel Science Foundation under ISF Grant 1050/16 and the Knut and Alice Wallenberg foundation. 

%	{\small
%		\bibliographystyle{myamsalpha}
%		\bibliography{ref}}
%
	
	{\small
		\bibliographystyle{myamsalpha}
		\bibliography{ref}}

\end{document}